\documentclass[12pt]{amsart}%
\usepackage{graphicx}
\usepackage{amscd}
\usepackage{amsmath}
\usepackage{amsfonts}
\usepackage{amssymb}%
\providecommand{\U}[1]{\protect\rule{.1in}{.1in}}
\theoremstyle{plain}

\theoremstyle{definition}

\newtheorem{remark}{Remark}[section]

\numberwithin{equation}{section}
\numberwithin{theorem}{section}
\begin{document}
\title[Fourth-order type equations]{The fourth-order type linear\\ordinary differential equations}
\author{W.N. Everitt}
\address{W.N. Everitt, School of Mathematics and Statistics, University of Birmingham,
Edgbaston, Birmingham B15 2TT, England, UK}
\email{w.n.everitt@bham.ac.uk}
\author{D. J. Smith}
\address{D. J. Smith, School of Mathematics and Statistics, University of Birmingham,
Edgbaston, Birmingham B15 2TT, England, UK}
\email{smithd@for.mat.bham.ac.uk}
\author{M. van Hoeij}
\address{M. van Hoeij, Department of Mathematics, Florida State University, 208 James
J. Love Building, Tallahassee, FL 32306-4510, USA}
\email{hoeij@math.fsu.edu}
\date{21 March 2006}
\subjclass[2000]{Primary: 33C15, 33D15, 33F10; Secondary: 33C05, 34B05.}
\keywords{Ordinary linear differential equations, special functions, Bessel-type,
Jacobi-type, Laguerre-type, Legendre-type equations.}

\begin{abstract}
This note reports on the recent advancements in the search for explicit
representation, in classical special functions, of the solutions of the
fourth-order linear ordinary differential equations named Bessel-type,
Jacobi-type, Laguerre-type, Legendre-type.

\end{abstract}
\maketitle
\tableofcontents

\section{Introduction\label{sec1}}

The four fourth-order type linear ordinary differential equations are named
after the four classical second-order equations of Bessel, Jacobi, Laguerre, Legendre.

The three fourth-order equations Jacobi-type, Laguerre-type, Legendre-type,
and their associated orthogonal polynomials, were first defined by H.L. Krall
in 1940, see \cite{HK} and \cite{HK1}, and later studied in detail by A.M.
Krall in 1981, see \cite{AK}, and by Koornwinder in 1984, see \cite{THK}.

The structured definition of the general-even order Bessel-type special
functions is dependent upon the Jacobi and Laguerre classical orthogonal
polynomials, and the Jacobi-type and Laguerre-type orthogonal polynomials; see
\cite{EM}.

The properties of the fourth-order Bessel-type functions have been studied in
the papers \cite{DEHLM}, \cite{EKLM} and \cite{EKLM1}.

This short paper reports on the recent progress that has been made in the
search for representation of solutions of these four type differential
equations in terms of classical special functions.

For the Bessel-type equation two independent solution were obtained, see
\cite{EM} and \cite{DEHLM}, and then a complete set of four independent
solutions in the papers \cite{MvH} and \cite{WNE1}; these solutions are
dependent upon the classical Bessel functions $J_{r},Y_{r},I_{r},K_{r}$ for
$r=0,1.$

Subsequently, using the techniques developed in \cite{MvH1}, solutions of the
three Jacobi-type, Laguerre-type, Legendre-type equations have been obtained
in terms of the solutions of the corresponding classical differential equations.

In all these four cases the verification of these representations as solutions
of the type differential equations, is a formal process using the
computational methods in the computer program Maple. Copies of the
corresponding Maple .mws files, for the Jacobi-type, Laguerre-type,
Legendre-type equations, are available; see Section \ref{sec6} below.

There are four subsequent sections devoted to reporting on the form of these
solutions, for each of the type fourth-order differential equations.

The four type differential equations are all written in the Lagrange symmetric
(formally self-adjoint) form with spectral parameter $\Lambda$ or $\lambda;$
if further properties of the solutions of these equations, involving the
theory of special functions, are to be studied, then these spectral notations
may have to be changed.

\section{Bessel-type differential equation\label{sec2}}

The fourth-order Bessel-type differential equation takes the form
\begin{equation}
(xy^{\prime\prime}(x))^{\prime\prime}-((9x^{-1}+8M^{-1}x)y^{\prime
}(x))^{\prime}=\Lambda xy(x)\ \text{for all}\ x\in(0,\infty) \label{eq2.1}%
\end{equation}
where $M\in(0,\infty)$ is a positive parameter and $\Lambda\in\mathbb{C},$ the
complex field, is a spectral parameter. The differential equation
(\ref{eq2.1}) is derived in the paper \cite[Section 1, (1.10a)]{EM}. The many
spectral properties of this equation in the weighted Hilbert function space
$L^{2}((0,\infty);x)$ are studied in the papers \cite{DEHLM}, \cite{EKLM} and
\cite{EKLM1}; for a collected account see the survey paper \cite{WNE1}.

Our knowledge of the special function solutions of the Bessel-type
differential equation (\ref{eq2.1}) is now more complete than at the time the
paper \cite{EM} was written. However, the results in \cite[Section 1,
(1.8a)]{EM}, with $\alpha=0,$ show that the function defined by%
\begin{equation}
J_{\lambda}^{0,M}(x):=[1+M(\lambda/2)^{2}]J_{0}(\lambda x)-2M(\lambda
/2)^{2}(\lambda x)^{-1}J_{1}(\lambda x)\ \text{for all}\ x\in(0,\infty),
\label{eq2.2}%
\end{equation}
is a solution of the differential equation (\ref{eq2.1}), for all $\lambda
\in\mathbb{C},$ and hence for all $\Lambda\in\mathbb{C},$ and all $M>0.$ Here:

\begin{itemize}
\item[$(i)$] the parameter $M>0$

\item[$(ii)$] the parameter $\lambda\in\mathbb{C}$

\item[$(iii)$] the spectral parameter $\Lambda$ and the parameter $M,$ in the
equation (\ref{eq2.1}), and the parameters $M$ and $\lambda,$ in the
definition (\ref{eq2.2}), are connected by the relationship%
\begin{equation}
\Lambda\equiv\Lambda(\lambda,M)=\lambda^{2}(\lambda^{2}+8M^{-1})\ \text{for
all}\ \lambda\in\mathbb{C}\ \text{and all}\ M>0 \label{eq2.3}%
\end{equation}

\item[$(iv)$] $J_{0}$ and $J_{1}$ are the classical Bessel functions (of the
first kind), see \cite[Chapter III]{GNW}.
\end{itemize}

Similar arguments to the methods given in \cite{EM} show that the function
defined by%
\begin{equation}
Y_{\lambda}^{0,M}(x):=[1+M(\lambda/2)^{2}]Y_{0}(\lambda x)-2M(\lambda
/2)^{2}(\lambda x)^{-1}Y_{1}(\lambda x)\ \text{for all}\ x\in(0,\infty),
\label{eq2.4}%
\end{equation}
is also an independent solution of the differential equation (\ref{eq2.1}),
for all $\lambda\in\mathbb{C},$ and hence for all $\Lambda\in\mathbb{C}$ and
all $M>0;$ here, again, $Y_{0}$ and $Y_{1}$ are classical Bessel functions (of
the second kind), see \cite[Chapter III]{GNW}.

These earlier studies of the fourth-order differential equation (\ref{eq2.1})
failed to find any explicit form of two linearly independent solutions,
additional to the solutions $J_{\lambda}^{0,M}$ and $Y_{\lambda}^{0,M}.$
However, results of van Hoeij, see \cite{MvH1} and \cite{MvH}, using the
computer algebra program Maple have yielded the required two additional
solutions, here given the notations of $I_{\lambda}^{0,M}$ and $K_{\lambda
}^{0,M},$ with explicit representation in terms of the classical modified
Bessel functions $I_{0},K_{0}$ and $I_{1},K_{1}.$ These two additional
solutions are defined as follows, where as far as possible we have followed
the notation used for the solutions $J_{\lambda}^{0,M}$ and $Y_{\lambda}%
^{0,M},$

\begin{itemize}
\item[$(i)$] given $\lambda\in\mathbb{C}$, with $\arg(\lambda)\in\lbrack
0,2\pi),$ $M\in(0,\infty)$ and using the principal value of $\sqrt{\cdot},$
define%
\begin{equation}
c\equiv c(\lambda,M):=\sqrt{\lambda^{2}+8M^{-1}}\ \text{and}\ d\equiv
d(\lambda,M):=1+M(\lambda/2)^{2} \label{eq2.5a}%
\end{equation}

\item[$(ii)$] define the solution, for all $x\in(0,\infty),$%
\begin{align}
I_{\lambda}^{0,M}(x)  &  :=-dI_{0}(cx)+\tfrac{1}{2}cMx^{-1}I_{1}%
(cx)\label{eq2.5b}\\
&  :=-[1+M(\lambda/2)^{2}]I_{0}\left(  x\sqrt{\lambda^{2}+8M^{-1}}\right)
+\sqrt{M\left(  2+M(\lambda/2)^{2}\right)  }x^{-1}I_{1}\left(  x\sqrt
{\lambda^{2}+8M^{-1}}\right)  \label{eq2.5d}%
\end{align}

\item[$(iii)$] define the solution, for all $x\in(0,\infty),$
\begin{align}
K_{\lambda}^{0,M}(x)  &  :=dK_{0}(cx)+\tfrac{1}{2}cMx^{-1}K_{1}%
(cx)\label{eq2.5c}\\
&  :=[1+M(\lambda/2)^{2}]K_{0}\left(  x\sqrt{\lambda^{2}+8M^{-1}}\right)
+\sqrt{M\left(  2+M(\lambda/2)^{2}\right)  }x^{-1}K_{1}\left(  x\sqrt
{\lambda^{2}+8M^{-1}}\right)  . \label{eq2.5e}%
\end{align}

\end{itemize}

\begin{remark}
\label{rem2.1}We have

\begin{enumerate}
\item The four linearly independent solutions $J_{\lambda}^{0,M},Y_{\lambda
}^{0,M},I_{\lambda}^{0,M},K_{\lambda}^{0,M}$ provide a basis for all solutions
of the original differential equation (\ref{eq2.1}), subject to the
$(\Lambda,\lambda)$ connection given in (\ref{eq2.3}).

\item These four solutions are real-valued on their domain $(0,\infty)$ for
all $\lambda\in\mathbb{R}.$

\item The domain $(0,\infty)$ of the solutions $J_{\lambda}^{0,M}$ and
$I_{\lambda}^{0,M}$ can be extended to the closed half-line $[0,\infty)$ with
the properties%
\[
J_{\lambda}^{0,M}(0)=I_{\lambda}^{0,M}(0)=1\ \text{for all}\ \lambda
\in\mathbb{R}\ \text{and all}\ M\in(0,\infty).
\]

\end{enumerate}
\end{remark}

\section{Laguerre-type differential equation\label{sec3}}

The Laguerre-type differential equation was discovered by H.L. Krall, see
\cite{HK} and \cite{HK1}, and subsequently studied by other authors, see
\cite{AK}, \cite{THK} and \cite{WNE}.

The differential equation may be written in two forms; here the parameter
$A\in(0,\infty)$ and the spectral parameter $\lambda\in\mathbb{C}$:

\begin{enumerate}
\item The Frobenius form:%
\begin{equation}%
\begin{array}
[c]{r}%
x^{2}y^{(4)}(x)+(-2x^{2}+4x)y^{(3)}(x)+(x^{2}-(2A+6)x)y^{\prime\prime}(x)\\
+((2A+2)x-2A)y^{\prime}(x)=\lambda y(x)\ \text{for all}\ x\in(0,\infty).
\end{array}
\label{eq3.1}%
\end{equation}

\item The Lagrange symmetric form (formally self-adjoint form):%
\begin{equation}%
\begin{array}
[c]{r}%
\left(  (x^{2}\exp(-x)y^{\prime\prime}(x)\right)  ^{\prime\prime}-\left(
((2A+2)x+2)\exp(-x)y^{\prime}(x)\right)  ^{\prime}\\
=\lambda\exp(-x)y(x)\ \text{for all}\ x\in(0,\infty).
\end{array}
\label{eq3.2}%
\end{equation}

\end{enumerate}

The spectral properties of the equation (\ref{eq3.2}) are considered in the
weighted Hilbert function space $L^{2}((0,\infty);\exp(\cdot)).$

The van Hoeij method of searching for factors and solutions of this
differential equation depends on the use of the computer algebra program
Maple; see \cite{MvH1} and \cite{MvH}.

Following the solutions of the fourth-order Bessel-type obtained in
\cite[Section 2]{EM}, and then in \cite{WNE1} on using the methods of
\cite{MvH}, it has proved possible to express all solutions of the equations
(\ref{eq3.1}) and (\ref{eq3.2}) in terms of the confluent hypergeometric functions.

The special functions involved are the Whittaker functions $M_{\kappa,\mu}(z)$
and $W_{\kappa,\mu}(z)$ as defined in the compendium \cite[Chapter 13,
Sections 13.1.31 to 13.1.34]{AS}.

Define the functions $\Phi$ and $\Psi$, for the range of variables given in
\cite[Chapter 13]{AS}:%
\begin{equation}
\Phi(\kappa,\mu;z):=M_{\kappa,\mu}(z) \label{eq3.2a}%
\end{equation}
\begin{equation}
\Psi(\kappa,\mu;z):=W_{\kappa,\mu}(z). \label{eq3.2b}%
\end{equation}

Given the parameters $A\in(0,\infty)$ and $\lambda\in\mathbb{C}$ define%
\begin{equation}
\Gamma(\lambda,A):=\sqrt{4A^{2}+4A+1+4\lambda}. \label{eq3.2c}%
\end{equation}

We designate the four linearly independent solutions of (\ref{eq3.1}),
equivalently (\ref{eq3.2}), again the parameters $A\in(0,\infty)$ and
$\lambda\in\mathbb{C},$ by
\[
L_{r}(\lambda,A;x)\ \text{for all}\ x\in(0,\infty)\ \text{and}\ r=1,2,3,4.
\]

These solutions are defined as follows, where $d/dx$ denotes differentiation
with respect to the independent variable $x,$

\begin{itemize}
\item[$L_{1}$]
\begin{equation}%
\begin{array}
[c]{cr}%
L_{1}(\lambda,A;x):= & \left(  \frac{1}{2}+\frac{1}{2}\Gamma(\lambda
,A)\right)  x^{-1/2}\exp(\frac{1}{2}x)\Phi(-A-\frac{1}{2}\Gamma(\lambda
,A),0;x)\\
& -\dfrac{d}{dx}\left(  x^{-1/2}\exp(\frac{1}{2}x)\Phi(-A-\frac{1}{2}%
\Gamma(\lambda,A),0;x)\right)  .
\end{array}
\label{eq3.3}%
\end{equation}

\item[$L_{2}$]
\begin{equation}%
\begin{array}
[c]{cr}%
L_{2}(\lambda,A;x):= & \left(  \frac{1}{2}-\frac{1}{2}\Gamma(\lambda
,A)\right)  x^{-1/2}\exp(\frac{1}{2}x)\Phi(-A+\frac{1}{2}\Gamma(\lambda
,A),0;x)\\
& -\dfrac{d}{dx}\left(  x^{-1/2}\exp(\frac{1}{2}x)\Phi(-A+\frac{1}{2}%
\Gamma(\lambda,A),0;x)\right)  .
\end{array}
\label{eq3.4}%
\end{equation}

\item[$L_{3}$]
\begin{equation}%
\begin{array}
[c]{cr}%
L_{3}(\lambda,A;x):= & \left(  \frac{1}{2}+\frac{1}{2}\Gamma(\lambda
,A)\right)  x^{-1/2}\exp(\frac{1}{2}x)\Psi(-A-\frac{1}{2}\Gamma(\lambda
,A),0;x)\\
& -\dfrac{d}{dx}\left(  x^{-1/2}\exp(\frac{1}{2}x)\Psi(-A-\frac{1}{2}%
\Gamma(\lambda,A),0;x)\right)  .
\end{array}
\label{eq3.5}%
\end{equation}

\item[$L_{4}$]
\begin{equation}%
\begin{array}
[c]{cr}%
L_{4}(\lambda,A;x):= & \left(  \frac{1}{2}-\frac{1}{2}\Gamma(\lambda
,A)\right)  x^{-1/2}\exp(\frac{1}{2}x)\Psi(-A+\frac{1}{2}\Gamma(\lambda
,A),0;x)\\
& -\dfrac{d}{dx}\left(  x^{-1/2}\exp(\frac{1}{2}x)\Psi(-A+\frac{1}{2}%
\Gamma(\lambda,A),0;x)\right)  .
\end{array}
\label{eq3.6}%
\end{equation}

\end{itemize}

\section{Legendre-type differential equation\label{sec4}}

The Legendre-type differential equation was discovered by H.L. Krall, see
\cite{HK} and \cite{HK1}, and subsequently studied by other authors, see
\cite{AK}, \cite{THK} and \cite{WNE}.

The differential equation may be written in two forms; here the parameter
$A\in(0,\infty)$ and the spectral parameter $\lambda\in\mathbb{C}$:

\begin{enumerate}
\item The Frobenius form:%
\begin{equation}%
\begin{array}
[c]{r}%
(x^{2}-1)^{2}y^{(4)}(x)+8x(x^{2}-1)y^{(3)}(x)+(4A+12)(x^{2}-1)y^{\prime\prime
}(x)\\
+8Axy^{\prime}(x)=\lambda y(x)\ \text{for all}\ x\in(-1,+1).
\end{array}
\label{eq4.1}%
\end{equation}

\item The Lagrange symmetric form (formally self-adjoint form):%
\begin{equation}%
\begin{array}
[c]{r}%
\left(  (1-x^{2})^{2}y^{\prime\prime}(x)\right)  ^{\prime\prime}-\left(
(8+4A(1-x^{2}))y^{\prime}(x)\right)  ^{\prime}\\
=\lambda y(x)\ \text{for all}\ x\in(-1,+1).
\end{array}
\label{eq4.2}%
\end{equation}

\end{enumerate}

The spectral properties of the equation (\ref{eq4.2}) are considered in the
Hilbert function space $L^{2}(-1,+1).$

The van Hoeij method of searching for factors and solutions of this
differential equation depends on the use of the computer algebra program
Maple; see \cite{MvH1} and \cite{MvH}.

Following the solutions of the fourth-order Bessel-type obtained in
\cite[Section 2]{EM}, and then in \cite{WNE1} on using the methods of
\cite{MvH}, it has proved possible to express all solutions of the equations
(\ref{eq4.1}) and (\ref{eq4.2}) in terms of the confluent hypergeometric functions.

The special functions involved are the Legendre functions $P_{\nu}^{\mu}(z)$
and $Q_{\nu}^{\mu}(z)$ as defined in the compendium \cite[Chapter 8, Section
8.1]{AS}.

Define the functions $P(\nu,x)$ and $Q(\nu,x)$, for the range of variables
given in \cite[Chapter 8]{AS}, by:%
\begin{equation}
P(\nu,x):=P_{\nu}^{0}(x) \label{eq4.3}%
\end{equation}
\begin{equation}
Q(\nu,x):=Q_{\nu}^{0}(z). \label{eq4.4}%
\end{equation}

Given the parameters $A\in(0,\infty)$ and $\lambda\in\mathbb{C}$ define,
noting the use of the $\pm$ symbol,%
\begin{equation}
\Gamma^{\pm}(\lambda,A):=\pm\sqrt{4A^{2}-4A+1+\lambda}\ \text{and}%
\ \Omega^{\pm}(\lambda,A):=\sqrt{5-8A+4\Gamma^{\pm}(\lambda,A)} \label{eq4.5}%
\end{equation}

We designate the four linearly independent solutions of (\ref{eq4.1}),
equivalently (\ref{eq4.2}), again the parameters $A\in(0,\infty)$ and
$\lambda\in\mathbb{C},$ by
\begin{equation}
Le_{r}(\lambda,A;x)\ \text{for all}\ x\in(-1,+1)\ \text{and}\ r=1,2,3,4.
\label{eq4.6}%
\end{equation}

These solutions are defined as follows, where $d/dx$ denotes differentiation
with respect to the independent variable $x,$%
\begin{equation}%
\begin{array}
[c]{c}%
Le_{1}(\lambda,A;x):=\\
-\frac{1}{2}(1+\Omega^{+}(\lambda,A))xP\left(  \frac{1}{2}\sqrt{9-8A+4\Gamma
^{+}(\lambda,A)+4\Omega^{+}(\lambda,A)}-\frac{1}{2},x\right)  +\\
\lbrack-(\lambda+3-4A+4A^{2})+\Omega^{+}(\lambda,A)-3\Gamma^{+}(\lambda
,A)+\Omega^{+}(\lambda,A)\Gamma^{+}(\lambda,A)+(\lambda+4A+4A^{2})x^{2}%
]\times\\
(\lambda+4A+4A^{2})^{-1}\dfrac{d}{dx}P(\frac{1}{2}\sqrt{9-8A+4\Gamma
^{+}(\lambda,A)+4\Omega^{+}(\lambda,A)}-\frac{1}{2},x).
\end{array}
\label{eq4.7}%
\end{equation}%
\begin{equation}%
\begin{array}
[c]{c}%
Le_{2}(\lambda,A;x):=\\
-\frac{1}{2}(1+\Omega^{-}(\lambda,A))xP\left(  \frac{1}{2}\sqrt{9-8A+4\Gamma
^{-}(\lambda,A)+4\Omega^{-}(\lambda,A)}-\frac{1}{2},x\right)  +\\
\lbrack-(\lambda+3-4A+4A^{2})+\Omega^{-}(\lambda,A)-3\Gamma^{-}(\lambda
,A)+\Omega^{-}(\lambda,A)\Gamma^{-}(\lambda,A)+(\lambda+4A+4A^{2})x^{2}%
]\times\\
(\lambda+4A+4A^{2})^{-1}\dfrac{d}{dx}P(\frac{1}{2}\sqrt{9-8A+4\Gamma
^{-}(\lambda,A)+4\Omega^{-}(\lambda,A)}-\frac{1}{2},x).
\end{array}
\label{eq4.8}%
\end{equation}%
\begin{equation}%
\begin{array}
[c]{c}%
Le_{3}(\lambda,A;x):=\\
-\frac{1}{2}(1+\Omega^{+}(\lambda,A))xQ\left(  \frac{1}{2}\sqrt{9-8A+4\Gamma
^{+}(\lambda,A)+4\Omega^{+}(\lambda,A)}-\frac{1}{2},x\right)  +\\
\lbrack-(\lambda+3-4A+4A^{2})+\Omega^{+}(\lambda,A)-3\Gamma^{+}(\lambda
,A)+\Omega^{+}(\lambda,A)\Gamma^{+}(\lambda,A)+(\lambda+4A+4A^{2})x^{2}%
]\times\\
(\lambda+4A+4A^{2})^{-1}\dfrac{d}{dx}Q(\frac{1}{2}\sqrt{9-8A+4\Gamma
^{+}(\lambda,A)+4\Omega^{+}(\lambda,A)}-\frac{1}{2},x).
\end{array}
\label{eq4.9}%
\end{equation}%
\begin{equation}%
\begin{array}
[c]{c}%
Le_{4}(\lambda,A;x):=\\
-\frac{1}{2}(1+\Omega^{-}(\lambda,A))xQ\left(  \frac{1}{2}\sqrt{9-8A+4\Gamma
^{-}(\lambda,A)+4\Omega^{-}(\lambda,A)}-\frac{1}{2},x\right)  +\\
\lbrack-(\lambda+3-4A+4A^{2})+\Omega^{-}(\lambda,A)-3\Gamma^{-}(\lambda
,A)+\Omega^{-}(\lambda,A)\Gamma^{-}(\lambda,A)+(\lambda+4A+4A^{2})x^{2}%
]\times\\
(\lambda+4A+4A^{2})^{-1}\dfrac{d}{dx}Q(\frac{1}{2}\sqrt{9-8A+4\Gamma
^{-}(\lambda,A)+4\Omega^{-}(\lambda,A)}-\frac{1}{2},x).
\end{array}
\label{eq4.10}%
\end{equation}

\begin{enumerate}
\item It has been verified, using Maple 9, that $Le_{1},Le_{2},Le_{3},Le_{4}$
are solutions of the differential equations (\ref{eq4.1}) and (\ref{eq4.2})
for the domain $x\in(-1,+1),$ and the parameter values $A\in(0,\infty)$ and
$\lambda\in\mathbb{C};$ this verification follows from the original results of
van Hoeij.

\item The factor $(\lambda+4A+4A^{2})^{-1}$ which appears in the derivative
term of all four of the solutions $Le_{r},$ for $r=1,2,3,4,$ presents an
apparent singularity in these solutions when $\lambda=-(4A+4A^{2});$ however
Maple 9 determines that the other factors%
\[
\lbrack-(\lambda+3-4A+4A^{2})+\Omega^{\pm}(\lambda,A)-3\Gamma^{\pm}%
(\lambda,A)+\Omega^{\pm}(\lambda,A)\Gamma^{\pm}(\lambda,A)+(\lambda
+4A+4A^{2})x^{2}]
\]
both have a zero at this value $\lambda=-(4A+4A^{2})$.
\end{enumerate}

\section{Jacobi-type differential equation\label{sec5}}

This note concerns information about the van Hoeij solutions to the
fourth-order Jacobi-type linear ordinary differential equation.

The Jacobi-type differential equation was discovered by H.L. Krall, see
\cite{HK} and \cite{HK1}, and subsequently studied by other authors, see
\cite{AK}, \cite{THK} and \cite{WNE}.

The differential equation may be written in two forms; here the parameters
$A\in(0,\infty)$ and $\alpha\in(-1,\infty),$ and the spectral parameter
$\lambda\in\mathbb{C}$:

\begin{enumerate}
\item The Frobenius form:%
\begin{equation}
\left\{
\begin{array}
[c]{l}%
(1-x^{2})^{2}y^{(4)}(x)-2(1-x^{2})((\alpha+4)x+\alpha)y^{(3)}(x)\\
+(1+x)\left(  (4A2^{\alpha}+\alpha^{2}+9\alpha+14)x+(-4A2^{\alpha}+\alpha
^{2}-3\alpha-10)\right)  y^{\prime\prime}(x)\\
+((4A\alpha2^{\alpha}+8A2^{\alpha}+2\alpha^{2}+6\alpha+4)x+(4A\alpha2^{\alpha
}+2\alpha^{2}+6\alpha+4)y^{\prime}(x)\\
\hspace{1in}\hspace{1in}=\lambda y(x)\ \text{for all}\ x\in(-1,+1).
\end{array}
\right.  \label{eq5.1}%
\end{equation}

\item The Lagrange symmetric form (formally self-adjoint form):%
\begin{equation}
\left\{
\begin{array}
[c]{l}%
J[y](x):=((1-x)^{\alpha+2}(1+x)^{2}y^{\prime\prime}(x))^{\prime\prime}\\
-((1-x)^{\alpha+1}((4A2^{\alpha}+2\alpha+2)x+4A2^{\alpha}+2\alpha+6)y^{\prime
}(x))^{\prime}\\
\hspace{1in}\hspace{1in}=\lambda(1-x)^{\alpha}y(x)\ \text{for all}%
\ x\in(-1,+1).
\end{array}
\right.  \label{eq5.2}%
\end{equation}

\end{enumerate}

The spectral properties of the equation (\ref{eq5.2}) are considered in the
weighted Hilbert function space $L^{2}((-1,+1);(1-x)^{\alpha}).$

The van Hoeij method of searching for factors and solutions of this
differential equation depends on the use of the computer algebra program
Maple; see \cite{MvH1} and \cite{MvH}.

Following the solutions of the fourth-order Bessel-type obtained in
\cite[Section 2]{EM}, and then in \cite{WNE1} on using the methods of
\cite{MvH}, it has proved possible to express all solutions of the equations
(\ref{eq5.1}) and (\ref{eq5.2}) in terms of hypergeometric functions.

The special functions involved are the JacobiP functions; these functions are
defined in terms of the hypergeometric functions $_{p}F_{q};$ for certain
values of the parameters involved in this definition, they reduce to the
classical Jacobi polynomials $P_{n}^{(\alpha,\beta)}.$ The notation JacobiP is
not to be found in \cite[Chapter 22]{AS} but is given in the program Maple9;
here we shorten the notation JacobiP and define the function $JP$ as follows,
in terms of the classical $\Gamma$ and $_{2}F_{1}$ functions:%
\begin{equation}
\left\{
\begin{array}
[c]{l}%
\text{JacobiP}(\nu,\alpha,\beta;z)\equiv JP(\nu,\alpha,\beta;z)\\
\hspace{1.1in}\;\;\;\;:=\dfrac{\Gamma(\nu+\alpha+1)}{\Gamma(\nu+1)\Gamma
(\alpha+1)}{}_{2}F_{1}(-\nu,\nu+\alpha+\beta+1;\alpha+1;(1-z)/2).
\end{array}
\right.  \label{eq5.3}%
\end{equation}
Here, the independent variables $\nu,\alpha,\beta,z$ are in the complex field
$\mathbb{C},$ but for the results in this note we make the restrictions%
\[
\nu\in\mathbb{C},\;\alpha\in(-1,\infty),\;\beta=0,\;z\in\mathbb{R}.
\]

When $\nu=n\in\mathbb{N}_{0}$ then we have the connection, see \cite[Section
22.2]{AS} and \cite[Lecture 2, (2.2)]{RA},%
\[
P_{n}^{(\alpha,\beta)}(x)=JP(n,\alpha,\beta;x)\ \text{for all}\ \alpha
,\beta\in(-1,\infty)\ \text{and all}\ x\in(-1,+1).
\]

The Jacobi-type orthogonal polynomials were discovered by H.L. Krall, see
\cite{HK} and \cite{HK1}, and later studied by A.M. Krall, T.H. Koornwinder
and L.L. Littlejohn. These polynomials are denoted by $P_{n,A}^{\alpha}(x)$
for $\alpha\in(-1,\infty),A\in(0,\infty),n\in\mathbb{N}_{0}$ and
$x\in(-1,+1);$ each $P_{n,A}^{\alpha}$ is real-valued on $(-1,+1)$ and of
degree $n.$

These polynomials have the properties:

\begin{itemize}
\item[$(i)$] For all $n\in\mathbb{N}_{0}$ the polynomial $P_{n,A}^{\alpha}$ is
a solution of the differential equation, see (\ref{eq5.2})%
\begin{equation}
J[P_{n,A}^{\alpha}](x)=\lambda_{n}(\alpha,A)(1-x)^{\alpha}P_{n,A}^{\alpha
}(x)\ \text{for all}\ x\in(-1,+1) \label{eq5.4}%
\end{equation}
where the eigenvalue $\lambda_{n}(\alpha,A)$ is determined by%
\begin{equation}
\lambda_{n}(\alpha,A)=n(n+\alpha+1)(n^{2}+(\alpha+1)n+4A2^{\alpha}%
+\alpha)\ \text{for all}\ n\in\mathbb{N}_{0}. \label{eq5.5}%
\end{equation}

\item[$(ii)$] The collection $\{P_{n,A}^{\alpha}(\cdot):n\in\mathbb{N}_{0}\}$
is a complete, orthogonal set in the Lebesgue-Stieltjes function space
$L^{2}\left(  (-1,+1);\mu_{\text{Jac}}\right)  ,$ where the Borel measure
$\mu_{\text{Jac}}$ is determined by the monotonic non-decreasing function
$\hat{\mu}_{\text{Jac}}$, for all $\alpha\in(-1,\infty)$ and all
$A\in(0,\infty),$%
\begin{equation}
\left\{
\begin{array}
[c]{ccl}%
\hat{\mu}_{\text{Jac}}(x) & := & -1/2\ \text{for all}\ x\in(-\infty,-1)\\
& := & \dfrac{A}{2(\alpha+1)}\left[  2^{\alpha+1}-(1-x)^{\alpha+1}\right]
\ \text{for all}\ x\in\lbrack-1,+1]\\
& := & \dfrac{A2^{\alpha+1}}{2(\alpha+1)}\ \text{for all}\ x\in(+1,\infty).
\end{array}
\right.  \label{eq5.6}%
\end{equation}

\end{itemize}

The orthogonality of the set $\{P_{n,A}^{\alpha}:n\in\mathbb{N}_{0}\}$ then
appears as%
\begin{equation}%
\begin{array}
[c]{l}%
{\displaystyle\int_{\lbrack-1,+1]}}
P_{m,A}^{\alpha}(x)P_{n,A}^{\alpha}(x)~d\mu_{\text{Jac}}(x)\\
\hspace{0.5in}=\dfrac{1}{2}P_{m,A}^{\alpha}(-1)P_{n,A}^{\alpha}(-1)+\dfrac
{1}{2}A%
{\displaystyle\int_{-1}^{+1}}
P_{m,A}^{\alpha}(x)P_{n,A}^{\alpha}(x)(1-x)^{\alpha}~dx\\
\hspace{1in}=0\ \text{for all}\ m,n\in\mathbb{N}_{0}\ \text{with}\ m\neq n.
\end{array}
\label{eq5.7}%
\end{equation}

The van Hoeij solutions for the general differential equation (\ref{eq5.2})
begin with solutions to the set of special equations (\ref{eq5.4}) when the
general spectral parameter $\lambda$ takes the value of one of the eigenvalues
$\{\lambda_{n}(\alpha,A):n\in\mathbb{N}_{0}\}.$

Given $n\in\mathbb{N}_{0}$ define, for all $x\in(-1,+1),$ with $JP$ given by
(\ref{eq5.3}):%
\begin{equation}%
\begin{array}
[t]{ccc}%
S_{1,n}(x) & := & (n\alpha+2A2^{\alpha}+n+n^{2})JP(n,\alpha,0;x)\\
&  & +(1-x)\dfrac{d}{dx}JP(n,\alpha,0;x).
\end{array}
\label{eq5.8}%
\end{equation}
Then the van Hoeij methods determine that $S_{1,n}$ is a solution of
(\ref{eq5.1}) for this value of $n,$ \textit{i.e.} for $\lambda=\lambda
_{n}(\alpha,A).$

\begin{remark}
\label{rem4.1}From the definition $(\ref{eq5.8})$ it would appear the this
solution $S_{1,n}$ is linearly dependent upon the Jacobi-type polynomial
$P_{n,A}^{\alpha}.$
\end{remark}

To proceed further, given $\lambda\in\mathbb{C},$ let $\{\rho_{r}(\lambda
)\in\mathbb{C}:r=1,2,3,4\}$ denote the roots of the quartic polynomial%
\begin{equation}
\rho(\rho+\alpha+1)(\rho^{2}+(\alpha+1)\rho+4A2^{\alpha}+\alpha)-\lambda=0.
\label{eq5.9}%
\end{equation}
Maple gives explicit representations for these roots as follows:%
\begin{equation}
\left\{
\begin{array}
[c]{c}%
\rho_{1}(\alpha,A,\lambda)=\frac{1}{2}\left(  -\alpha-1+\sqrt{\alpha
^{2}+1-8A2^{\alpha}+2\xi(\alpha,A,\lambda)}\right) \\
\rho_{2}(\alpha,A,\lambda)=\frac{1}{2}\left(  -\alpha-1-\sqrt{\alpha
^{2}+1-8A2^{\alpha}+2\xi(\alpha,A,\lambda)}\right) \\
\rho_{3}(\alpha,A,\lambda)=\frac{1}{2}\left(  -\alpha-1+\sqrt{\alpha
^{2}+1-8A2^{\alpha}-2\xi(\alpha,A,\lambda)}\right) \\
\rho_{4}(\alpha,A,\lambda)=\frac{1}{2}\left(  -\alpha-1-\sqrt{\alpha
^{2}+1-8A2^{\alpha}-2\xi(\alpha,A,\lambda)}\right)  ,
\end{array}
\right.  \label{eq5.10}%
\end{equation}
where%
\[
\xi(\alpha,A,\lambda):=\sqrt{\alpha^{2}+8\alpha A2^{\alpha}+16A^{2}2^{2\alpha
}+4\lambda}.
\]

The van Hoeij solutions to the Jacobi-type equation, see (\ref{eq5.10}),%
\begin{equation}%
\begin{array}
[c]{c}%
J[y](x)=\lambda(1-x)^{\alpha}y(x)\ \text{for all}\ x\in(-1,+1)
\end{array}
\label{eq5.11}%
\end{equation}
now develop as follows.

With $\lambda\in\mathbb{C},$ replace $n$ in the solution $S_{1,n}$ with the
first root $\rho_{1}(\lambda)$ taken from (\ref{eq5.10}); this yields a
solution $\mathcal{J}_{1}$ in terms of the function $JP,$ for all $\alpha
\in(-1,\infty)$ and all $A\in(0,\infty),$%
\begin{equation}
\left\{
\begin{array}
[c]{lll}%
\mathcal{J}_{1}(\alpha,A,\lambda;x) & := & \dfrac{1}{2}\left[  \alpha
-\xi(\alpha,A,\lambda)\right]  JP(\rho_{1}(\alpha,A,\lambda),\alpha,0;x)\\
&  & \hspace{0.25in}-(1-x)\dfrac{d}{dx}JP(\rho_{1}(\alpha,A,\lambda
),\alpha,0;x)\ \text{for all}\ x\in(-1,+1).
\end{array}
\right.  \label{eq5.12}%
\end{equation}

Similarly, for $\lambda\in\mathbb{C},$ there is a solution $\mathcal{J}_{2}$
derived from the third root $\rho_{3}(\lambda),$ for all $\alpha\in
(-1,\infty)$ and all $A\in(0,\infty),$%
\begin{equation}
\left\{
\begin{array}
[c]{lll}%
\mathcal{J}_{2}(\alpha,A,\lambda;x) & := & \dfrac{1}{2}\left[  \alpha
+\xi(\alpha,A,\lambda)\right]  JP(\rho_{3}(\alpha,A,\lambda),\alpha,0;x)\\
&  & \hspace{0.25in}-(1-x)\dfrac{d}{dx}JP(\rho_{3}(\alpha,A,\lambda
),\alpha,0;x)\ \text{for all}\ x\in(-1,+1).
\end{array}
\right.  \label{eq5.14}%
\end{equation}

Given $n\in\mathbb{N}_{0}$ define, for all $x\in(-1,+1),$ there is a second
solution $S_{2,n},$ similar to $S_{1,n}$ of (\ref{eq4.1}), defined by%
\begin{equation}%
\begin{array}
[c]{ccc}%
S_{2,n}(x) & := & \left[  \left(  (n+1)\alpha+A2^{\alpha+1}+n+n^{2}\right)
\right]  (x-1)^{-\alpha}JP(-1-n,-\alpha,0;x)\\
&  & -(x-1)^{1-\alpha}\dfrac{d}{dx}JP(-1-n,-\alpha,0;x).
\end{array}
\label{eq5.16}%
\end{equation}
Then the van Hoeij methods determine that $S_{2,n}$ is a solution of
(\ref{eq5.1}) for this value of $n,$ \textit{i.e.} for $\lambda=\lambda
_{n}(\alpha,A).$

With $\lambda\in\mathbb{C},$ replace $n$ in the solution $S_{2,n}$ with the
first root $\rho_{1}(\lambda)$ taken from (\ref{eq5.10}); this yields a
solution $\mathcal{J}_{3}$ in terms of the function $JP,$ for all $\alpha
\in(-1,\infty)$ and all $A\in(0,\infty),$%
\begin{equation}
\left\{
\begin{array}
[c]{ccl}%
\mathcal{J}_{3}(a,A,\lambda;x) & := & -\dfrac{1}{2}[\alpha+\xi(\alpha
,A,\lambda)](-1+x)^{-\alpha}JP(-1-\rho_{1}(\alpha,A,\lambda),-\alpha,0;x)\\
&  & \hspace{0.25in}+(-1+x)^{1-\alpha}\dfrac{d}{dx}JP(-1-\rho_{1}%
(\alpha,A,\lambda),-\alpha,0;x)\ \text{for all}\ x\in(-1,+1).
\end{array}
\right.  \label{eq5.17}%
\end{equation}

Similarly, for $\lambda\in\mathbb{C},$ there is a solution $\mathcal{J}_{4}$
derived from the third root $\rho_{3}(\lambda),$ for all $\alpha\in
(-1,\infty)$ and all $A\in(0,\infty),$%
\begin{equation}
\left\{
\begin{array}
[c]{lll}%
\mathcal{J}_{4}(\alpha,A,\lambda;x) & := & -\dfrac{1}{2}[\alpha-\xi
(\alpha,A,\lambda)](-1+x)^{-\alpha}JP(-1-\rho_{3}(\alpha,A,\lambda
),-\alpha,0;x)\\
&  & \hspace{0.25in}+(-1+x)^{1-\alpha}\dfrac{d}{dx}JP(-1-\rho_{3}%
(\alpha,A,\lambda),-\alpha,0;x)\ \text{for all}\ x\in(-1,+1).
\end{array}
\right.  \label{eq5.19}%
\end{equation}

\section{Maple files\label{sec6}}

The formal validity of these solution results were originally obtained by van
Hoeij, following the methods given in \cite{MvH1} and \cite{MvH}, using the
computer program Maple. These results have been subsequently confirmed by
David Smith in the Maple program at the University of Birmingham.

The Maple .mws files for the following three cases are named as follows:

\begin{center}
\textbf{Jacobi-type-solutions.mws}

\textbf{Laguerre-type-solutions.mws}

\textbf{Legendre-type-solutions.mws}
\end{center}

Copies of these files may be obtained on application to David Smith at

\begin{center}
smithd@for.mat.bham.ac.uk
\end{center}

\section{Acknowledgements\label{sec7}}

The authors indicate here that the results reported on in this note follow
from the combined efforts of the many collaborators whose names are given in
the references below.

The author Norrie Everitt is particularly grateful for the advice and help
received from Dr David Smith, School of Mathematics at the University of
Birmingham, on the application and use of the computer program Maple.

The author Mark van Hoeij is supported by National Science Foundation grant 0511544.

\end{document}